\documentclass{ifacconf}

\usepackage{graphicx}      

\makeatletter
\let\old@ssect\@ssect 
\makeatother

\usepackage{natbib}
\usepackage{hyperref}

\makeatletter
\def\@ssect#1#2#3#4#5#6{%
  \NR@gettitle{#6}
  \old@ssect{#1}{#2}{#3}{#4}{#5}{#6}
}
\makeatother

\usepackage{amsmath,amsfonts,amssymb}
\usepackage{color}
\usepackage{caption}
\usepackage{subcaption}

\newcommand{\f}{{\mathrm{f}}}

\newcommand{\vd}{{\mathbf d}}
\newcommand{\vf}{{\mathbf f}}
\newcommand{\vg}{{\mathbf g}}
\newcommand{\vm}{{\mathbf m}}
\newcommand{\vu}{{\mathbf u}}

\newcommand{\vw}{{\mathbf w}}
\newcommand{\vx}{{\mathbf x}}
\newcommand{\vy}{{\mathbf y}}
\newcommand{\vz}{{\mathbf z}}

\newcommand{\vF}{{\mathbf F}}

\newcommand{\vS}{{\mathbf S}}
\newcommand{\vU}{{\mathbf U}}
\newcommand{\vX}{{\mathbf X}}
\newcommand{\vY}{{\mathbf Y}}
\newcommand{\vZ}{{\mathbf Z}}

\newcommand{\vmu}{{\boldsymbol{\mu}}}
\newcommand{\vnu}{{\boldsymbol{\nu}}}
\newcommand{\vSigma}{{\boldsymbol{\Sigma}}}
\newcommand{\vupsilon}{{\boldsymbol{\upsilon}}}

\newcommand{\Eb}{\mathbb{E}}

\newcommand{\Rb}{\mathbb{R}}
\newcommand{\Sb}{\mathbb{S}}

\usepackage{calrsfs}
\DeclareMathAlphabet{\pazocal}{OMS}{zplm}{m}{n}

\newcommand{\calD}{\pazocal{D}}

\newcommand{\calL}{\pazocal{L}}

\newcommand{\calN}{\pazocal{N}}

\newcommand{\calO}{\pazocal{O}}

\newcommand{\calU}{\pazocal{U}}

\DeclareMathOperator*{\argmin}{arg\,min}

\newtheorem{problem}{Problem}

\begin{document}
\begin{frontmatter}

\title{Nonlinear Covariance Steering using Variational Gaussian Process Predictive Models} 

\thanks[footnoteinfo]{This research has been supported in part by NSF awards ECCS-1924790, CCMI-1753687 and CMMI-1937957. Alexandros Tsolovikos acknowledges 
support by the A. Onassis Foundation scholarship.}
\thanks[git]{The code for this work is available at \href{https://github.com/alextsolovikos/greedyGPCS}{github.com/alextsolovikos/greedyGPCS}}

\author[UT]{Alexandros Tsolovikos} 
\author[UT]{Efstathios Bakolas} 

\address[UT]{The University of Texas at Austin, Austin, Texas 78712-1221 (\href{mailto:tsolovikos@utexas.edu}{tsolovikos@utexas.edu}, \href{mailto:bakolas@austin.utexas.edu}{bakolas@austin.utexas.edu})}


\begin{abstract}                
In this work, we consider the problem of steering the first two moments of the uncertain state of an \textit{unknown} discrete-time stochastic nonlinear system to a given terminal distribution in finite time. Toward that goal, first, a non-parametric predictive model is learned from a set of available training data points using stochastic variational Gaussian process regression: a powerful and scalable machine learning tool for learning distributions over arbitrary nonlinear functions. Second, we formulate a tractable nonlinear covariance steering algorithm that utilizes the Gaussian process predictive model to compute a feedback policy that will drive the distribution of the state of the system close to the goal distribution. In particular, we implement a greedy covariance steering control policy that linearizes at each time step the Gaussian process model around the latest predicted mean and covariance, solves the linear covariance steering control problem, and applies only the first control law. The state uncertainty under the latest feedback control policy is then propagated using the unscented transform with the learned Gaussian process predictive model and the algorithm proceeds to the next time step. Numerical simulations illustrating the main ideas of this paper are also presented.

\end{abstract}

\begin{keyword}
nonparametric methods, nonlinear system identification, stochastic system identification, covariance steering, stochastic optimal control problems
\end{keyword}

\end{frontmatter}

\section{Introduction}
In this paper, we consider the finite-horizon covariance steering problem for discrete-time stochastic nonlinear systems described by \textit{non-parametric Gaussian process} models. In particular, we consider the problem of learning sparse stochastic variational Gaussian process (SVGP) predictive models for stochastic nonlinear systems from training data and then using the SVGP models for computing feedback control policies that steer the mean and covariance of the uncertain state of the underlying system to desired quantities at a given (finite) terminal time. This problem will be referred to as the \textit{Gaussian process-based nonlinear covariance steering} problem.

\textit{Literature Review:} 
Gaussian Processes (GP) [\cite{rasmussen2003gaussian}] are non-parametric regression models that describe distributions over functions and are ideal for learning predictive models for arbitrary nonlinear stochastic systems due to their flexibility and inherent ability to provide uncertainty estimates that capture both model uncertainties and process noise. GP regression models have been used extensively for learning predictive state models for dynamical systems [\cite{grimes2006dynamic, ko2007gaussian}] and observation models for state estimation [\cite{ko2007gp, ko2009gp}], as well as trajectory optimization [\cite{pan2014probabilistic, pan2015data}] and motion planning [\cite{mukadam2016gaussian}]. Inference using GP models is inherently dependent on the training data and the cost of inference with exact GPs scales with the cube of the number of training points. For that reason, a number of sparse approximations of GPs have been proposed in the literature, the most common using a set of ``inducing variables'' [\cite{quinonero2005unifying, titsias2009variational}]. Further scalability can be achieved by using stochastic variational inference [\cite{hoffman2013stochastic}], leading to sparse GP models that can be trained on large datasets with thousands or millions of data points using stochastic gradient descent, while retaining a small inference cost [\cite{hensman2013gaussian}].

The infinite-horizon covariance steering (or covariance control) problem for both continuous-time and discrete-time linear Gaussian systems has been studied extensively [\cite{p:skeltonIJC, p:skeltonTAC, p:Grig97}], while the finite-horizon problem has been addressed in \cite{p:georgiou15A,p:georgiou15B} for the continuous-time and in \cite{p:bakcdc16,p:PT2017} for the discrete-time case. Covariance control problems for incomplete and imperfect state information have also been studied in \cite{p:EBACC17, p:bakTAC2019}. Nonlinear density steering problems for feedback linearizable nonlinear systems were recently studied in \cite{p:caluya2019}, while an iterative covariance steering algorithm for nonlinear systems based on a linearization of the system along reference trajectories was presented in \cite{p:ridderhof2019}. Stochastic nonlinear model predictive control with probabilistic constraints can also be found in \cite{p:mesbah2014,p:SEHR2017}.

\textit{Main Contribution:} In this work, non-parametric state predictive models of discrete-time stochastic nonlinear systems with \textit{unknown dynamics} are learned using \textit{stochastic variational GP regression} and subsequently are used to control the mean and covariance of the state of the unknown systems in a greedy nonlinear covariance steering algorithm.


First, we introduce stochastic variational GP regression and present the process of learning SVGP predictive models for discrete-time dynamics from a set of training samples obtained by measuring the underlying stochastic nonlinear system of interest. Then, the non-parametric predictive model is used in a greedy finite-horizon covariance steering algorithm similar to the one presented in \cite{bakolas2020greedy}. 

The Gaussian process-based greedy nonlinear covariance steering algorithm consists of three steps. First, the SVGP model is linearized around the latest mean state prediction or estimation. Second, the feedback control policy that solves the linear Gaussian covariance steering problem from the current mean and covariance estimates to the target ones under the linearized system is computed, but only the first control law is executed.
Then, the state mean and covariance of the closed-loop system that results by applying the feedback control policy computed at the previous step are propagated to the next time step using the unscented transform [\cite{p:julier2002,p:julier2004}], modified to take into account the uncertainty estimates provided by the GP predictive model [\cite{ko2007gp}]. This three-step process is repeated in a shrinking-horizon model predictive control fashion until the final time step, when the terminal state mean and covariance should sufficiently approximate the goal quantities.

\textit{Structure of the paper:} The rest of the paper is organized as follows. In Section \ref{section: sysid}, the process of learning a predictive model from sample data points using SVGP regression is presented. The greedy nonlinear covariance steering problem for non-parametric GP predictive models is formulated in Section \ref{section: greedy cs}. Section \ref{section: results} presents numerical simulations and comparisons of the GP model with the analytic one. We conclude with remarks and directions for future research in Section \ref{section: conclusions}.

\textit{Notation:}
Given a random vector $\vx$, $\mathbb{E}[\vx]$ denotes its expected value (mean) and $\mathrm{Cov}[\vx]$ its covariance, where $\mathrm{Cov}[\vx]:=
\mathbb{E}[ (\vx - \mathbb{E}[\vx])(\vx - \mathbb{E}[\vx])^\top ]$. The space of real symmetric $n \times n$ matrices will be denoted by
$\mathbb{S}_n$. Furthermore, the convex cone of
$n\times n$ (symmetric) positive semi-definite and (symmetric)
positive definite matrices will be denoted by $\mathbb{S}^{+}_n$ and $\mathbb{S}^{++}_n$, respectively. 
Finite-length sequences are denoted as $\{\vx_1,\dots,\vx_N\} = \{\vx_i\}_{i=1}^N$.
The $i$-th element of a vector $\vx$ is denoted by $[\vx]_i$. Similarly, the $i$-th element of the $j$-th column of a matrix $M$ is denoted as $[M]_{ij}$. For a scalar-valued function $f(\cdot):\Rb^n \rightarrow \Rb$ and a sequence of vectors $\vX = \{\vx_i\}_{i=1}^N$, we define $f(\vX)$ as the vector with $[f(\vX)]_{i} = f(\vx_i)$. 
Similarly, if $k(\cdot,\cdot):\Rb^n \times \Rb^n \rightarrow \Rb$, then $k(\vX,\vX)$ is the matrix with elements $[k(\vX,\vX)]_{ij} = k(\vx_i,\vx_j)$.
Finally, if $\vx\in \mathbb{R}^n$ and $\vy\in \mathbb{R}^m$, then $[\vx;\vy] = [\vx^\top,\vy^\top]^\top \in \mathbb{R}^{n+m}$ will denote the vertical concatenation of $\vx$ and $\vy$.

\section{Stochastic Variational Gaussian Processes for Discrete-Time Dynamics}
\label{section: sysid}

\subsection{Sparse Variational Gaussian Process Regression}
\label{subsection: svgp}
Consider the vector $\vy$, where $y_i$ is a noisy observation of an unknown scalar-valued function $f(\cdot):\Rb^n \rightarrow \Rb$ at a known location $\vx_i$, for all $\vX = \{\vx_i\}_{i=1}^N$, and the measurement likelihood $p(y_i \mid f(\vx_i))$ is known. Let $\vf$ be the (unknown) vector containing the values of $f(\cdot)$ at the points $\vX$. We introduce a Gaussian prior on $f(\cdot)$, i.e. $f(\vx) \sim \calN\left(f(\vx) \mid m(\vx), k(\vx,\vx)\right)$, where $m(\cdot):\Rb^n \rightarrow \Rb$ is a chosen \textit{mean function} (e.g. zero, constant, or linear) and $k(\cdot,\cdot): \Rb^n\times \Rb^n \rightarrow \Rb$ is the \textit{kernel function} that measures the closeness between two input points and specifies the smoothness and continuity properties of the underlying function $f(\cdot)$. Now, the \textit{prior} over the vector $\vf$ can be written as
\begin{align}
    p(\vf ; \vX) = \calN\left(\vf \mid m(\vX), k(\vX,\vX) \right),
    \label{eq: f prior}
\end{align}
where the mean vector is defined as $[m(\vX)]_i = m(\vx_i)$ and the covariance is $[k(\vX,\vX)]_{ij} = k(\vx_i,\vx_j)$.

\subsubsection{Exact GP Inference:} The joint density of $\vy$ and $\vf$ is
\begin{align}
    p(\vy, \vf ; \vX) = p(\vy \mid \vf ; \vX) p(\vf ; \vX).
\end{align}
If the likelihood function is chosen to be Gaussian, e.g. $p(\vy \mid \vf ; \vX) = \calN\left(\vy \mid \vf, \sigma_{\epsilon}^2 I \right)$, then, the marginal likelihood 
\begin{align}
    p(\vy; \vX) &= \int p(\vy \mid \vf ; \vX) p(\vf ; \vX) d\vf \nonumber \\
    &= \calN\left(\vy \mid m(\vX), k(\vX,\vX) + \sigma_{\epsilon}^2 I \right)
\end{align}
is analytically computed and the \textit{hyperparameters} $\mathbf{\Theta} = \{\theta_m, \theta_k, \sigma_{\epsilon} \}$ that define the Gaussian process mean, kernel, and likelihood functions can be directly optimized by minimizing the negative log-likelihood of the training data, i.e.
\begin{align}
    \mathbf{\Theta}_{\mathrm{opt}} = \argmin_{\mathbf{\Theta}}  \left(-\mathrm{log}\ p(\vy; \vX)\right).
\end{align}

Prediction of $y_*$ on a new location $\vx_*$ is done by conditioning on the training data,
\begin{align}
    p(y_* ; \vx_*, \vy, \vX) &= \int p(y_*, \vy ; \vx_*, \vX) d\vy \nonumber \\
    &= \calN\left(y_* \mid \mu_*, \sigma_* \right),
    \label{eq: exact posterior}
\end{align}
where
\begin{align*}
    \mu_* &= m(\vx_*) + k(\vx_*, \vX) \left[ k(\vX,\vX) + \sigma_{\epsilon}^2 I \right]^{-1} \left(\vy - m(\vX)\right)\\
    \sigma_* &= k(\vx_*, \vx_*) - k(\vx_*, \vX) \left[ k(\vX,\vX) + \sigma_{\epsilon}^2 I \right]^{-1} k(\vX, \vx_*).
\end{align*}

The computational cost of inference is $\calO(N^3)$, which can be expensive when the number of training points $N$ is large.

\subsubsection{Sparse Variational GP Inference:} In order to reduce the inference cost of a GP model, we can use sparse approximations of Gaussian processes. Define a set of $M$ inducing locations $\vZ = \{\vz_i\}_{i=1}^M$, with $M \ll N$, where $M$ and $\vZ$ are parameters to be chosen. In addition, define the vector $\vu$ as $[\vu]_i = f(\vz_i)$. The joint density of $\vy$, $\vf$, and $\vu$ is
\begin{align}
    p(\vy, \vf, \vu) = p(\vy \mid \vf ; \vX) p(\vf \mid \vu ; \vX, \vZ) p(\vu; \vZ),
\end{align}
where $p(\vu; \vZ) = \calN\left(\vu \mid m(\vZ), k(\vZ,\vZ) \right)$ is the Gaussian prior on $\vu$ (similar to \eqref{eq: f prior}) and $p(\vf \mid \vu ; \vX, \vZ) = \calN\left(\vf \mid \tilde{\vmu}, \tilde{\vSigma} \right)$, with
\begin{align*}
    [\tilde{\vmu}]_i &= m(\vx_i) + k(\vx_i, \vZ) k(\vZ,\vZ)^{-1} \left(\vu - m(\vZ)\right),\\
    [\tilde{\vSigma}]_{ij} &= k(\vx_i, \vx_j) - k(\vx_i, \vZ) k(\vZ,\vZ)^{-1} k(\vZ, \vx_j).
\end{align*}
However, $\vu$ \textit{is unknown}, since $f(\cdot)$ is also unknown. Following \cite{hensman2013gaussian}, we choose a \textit{variational} posterior
\begin{align}
    q(\vf,\vu) = p(\vf \mid \vu ; \vX, \vZ) q(\vu),
    \label{eq: variational posterior}
\end{align}
where $q(\vu) = \calN \left(\vu \mid \vm, \vS \right)$ and $\vm$, $\vS$ are the parameters defining the variational distribution (along with $\vZ$). Since both terms in \eqref{eq: variational posterior} are Gaussian, we can get rid of $\vu$ by marginalizing over it, that is,
\begin{align}
    q(\vf \mid \vm, \vS ; \vX, \vZ) &= \int p(\vf \mid \vu ; \vX, \vZ) q(\vu) d\vu \nonumber \\
    &= \calN\left(\vf \mid \vmu, \vSigma \right),
\end{align}
where, if we define the functions
\begin{align*}
    \mu_f(\vx_i) &:= m(\vx_i) + k(\vx_i, \vZ) \left[ k(\vZ,\vZ)\right]^{-1} \left(\vm - m(\vZ)\right),\\
    \vSigma_f(\vx_i,\vx_j) &:= k(\vx_i, \vx_j) \\
    &\hspace{-30pt}- k(\vx_i, \vZ) k(\vZ,\vZ)^{-1} \left[ k(\vZ,\vZ) - \vS \right] k(\vZ,\vZ)^{-1} k(\vZ, \vx_j),
\end{align*}
then $[\vmu]_i = \mu_f(\vx_i)$ and $[\vSigma]_{ij} = \vSigma_f(\vx_i,\vx_j)$. 

Once the variational parameters have been trained, \textit{predicting} the distribution of $y_*$ on a test location $\vx_*$ is simply
\begin{align*}
    p(y_* ; \vx_*, \vm, \vS, \vZ) = \calN\left( y_* \mid \mu_f(\vx_*), \vSigma_f(\vx_*,\vx_*) + \sigma_{\epsilon}^2\right).
\end{align*}
Now, only an $M \times M$ matrix needs to be inverted. The variational parameters ($\vZ$, $\vm$, and $\vS$), along with the hyperparameters $\mathbf{\Theta} = \{\theta_m, \theta_k, \sigma_{\epsilon} \}$, can be found by maximizing the lower bound $\calL$ on the marginal likelihood,
\begin{align}
    \log p(\vy \mid \vX) \geq \Eb_{q(\vf,\vu)} \left[ \log \frac{p(\vy,\vf,\vu)}{q(\vf,\vu)} \right] = \calL.
\end{align}
The lower bound can be factorized as
\begin{align}
    \calL = \sum_{i=1}^N \Eb_{q(f_i \mid \vm, \vS ; \vx_i, \vZ)} \left[ \log p(y_i \mid f_i) \right] - \mathrm{KL}\left[ q(\vu) \| p(\vu) \right],
    \label{eq: stochastic lower bound}
\end{align}
where KL denotes the Kullback-Leibler divergence. Note that the expectation can be computed analytically if the likelihood $p(y_i\mid f_i)$ is Gaussian. An immediate consequence of that choice is that, since the bound is the sum over the training data, we can perform stochastic inference through minibatch subsampling. This allows inference on large datasets and, more importantly, \textit{online learning} of the variational parameters.

\subsubsection{Multiple Outputs: } So far, the output $y_i \in \Rb$ has been a scalar. In the case of multiple outputs $\vy_i \in \Rb^D$, we can define the matrices $\vY$, $\vF$, and $\vU$ as the matrices containing the observation $\vy_i$ and function values $\vf(\vx_i)$ and $\vf(\vz_i)$ as their $i$-th rows. The latent functions are now $f_d(\cdot): \Rb^n \rightarrow \Rb$, for $d = 1,\dots,D$, and an independent sparse GP is learned for each function by maximizing a lower bound similar to \eqref{eq: stochastic lower bound}, but with $p(\vY, \vF, \vU) = \Pi_{d=1}^D p(\vy_d, \vf_d, \vu_d)$ and $q(\vF,\vU) = \Pi_{d=1}^D q(\vf_d,\vu_d)$ in place of $p(\vy,\vf,\vu)$ and $p(\vf,\vu)$, respectively.

\subsection{SVGP for Discrete-time Dynamics}
Consider a discrete-time dynamical system of the form
\begin{align}
    \vz_{t+1} = \vg(\vz_t, \vu_t) + \epsilon_t,
    \label{eq: exact dynamics}
\end{align}
where $\vz_t \in \Rb^{n_z}$ is the state at time step $t$, $\vu_t \in \Rb^{n_u}$ is the control input, and $\epsilon_t \in \Rb^{n_z}$ the i.i.d. additive Gaussian white noise, with $\vw_t \sim \calN\left(\epsilon_t \mid 0, \sigma_{\epsilon}^2 I \right)$.

Assume that the underlying dynamics $\vg(\cdot,\cdot)$ are unknown, but full-state measurements of the state transitions for given inputs are available for sampling (e.g. via experiments or simulation). In particular, assume that full-state \textit{observations} 
\begin{align}
    \vy_i = \vg(\vz_i, \vu_i) + \epsilon_i
    \label{eq: observations}
\end{align}
at known locations 
\begin{align}
    \vx_i = [\vz_i; \vu_i]
    \label{eq: inputs}
\end{align}
are available, that is, our dataset consists of $N$ triplets, $\calD = \{(\vy_i, \vz_i, \vu_i)\}_{i=1}^N$. Following Subsection \ref{subsection: svgp}, we can fit a multitask (multi-output) sparse variational GP to the measurements, in order to get a non-parametric approximate model of the dynamics. The observations and corresponding inputs to the SVGP are the ones defined in \eqref{eq: observations} and \eqref{eq: inputs}, respectively, the number of outputs is $D=n_z$ and the number of inputs (features) is $n = n_z + n_u$. Choosing an appropriate mean $m(\cdot)$ (e.g. zero, constant, or linear) and a kernel function $k(\cdot,\cdot)$ (typically, a squared exponential), the variational parameters and hyperparameters of the SVGP are learned by minimizing the negative of the lower bound, $-\calL$, via stochastic gradient descent on minibatches of $\calD$.

The learned transition SVGP model can now be defined as
\begin{align}
    \vz_{t+1} = G(\vz_t, \vu_t) + \vw_t,
    \label{eq: svgp model}
\end{align}
where
\begin{align}
    G(\vz_t, \vu_t) = \vmu_f([\vz_t; \vu_t])
    \label{eq: G definition}
\end{align}
is the mean of the next state and
\begin{align}
    \vw_t \sim \calN\left( \vw_t \mid 0, \Sigma_f([\vz_t; \vu_t], [\vz_t; \vu_t]) + \sigma_{\epsilon}^2 \right)
    \label{eq: w definition}
\end{align}
is the additive noise, the covariance of which captures not only the process noise, but also the model uncertainties.

\subsection{Linearization of the SVGP Dynamics}
Given a trained SVGP model like the one in \eqref{eq: svgp model}, if a linearization around a given state $\vz_*$ and input $\vu_*$ is necessary, it can be easily computed as 
\begin{align}
    \vz_{t+1} \approx A_* \vz_t + B_* \vu_t + \vd_*,
\end{align}
where
\begin{align}
    A_* &= \frac{\partial}{\partial \vz} G(\vz_*, \vu_*) = \frac{\partial}{\partial [\vz; \vu]} \vmu_f([\vz; \vu]) \begin{bmatrix} I_{n_z} \\ 0 \end{bmatrix}\bigg|_{ \substack{\vz = \vz_* \\ \vu = \vu_*}},\\
    B_* &= \frac{\partial}{\partial \vu} G(\vz_*, \vu_*) = \frac{\partial}{\partial [\vz; \vu]} \vmu_f([\vz; \vu]) \begin{bmatrix} 0 \\ I_{n_u} \end{bmatrix}\bigg|_{ \substack{\vz = \vz_* \\ \vu = \vu_*}},
    \label{eq: jacobians}
\end{align}
and
\begin{align*}
    \vd_* &= -A_* \vz_* - B_* \vu_* + G(\vz_*, \vu_*).
\end{align*}
For compactness, denote the linearization operation as
\begin{align}
    \{A_*, B_*, \vd_*\} = \mathrm{LIN}_{G}\{\vz_*, \vu_*\}.
\end{align}
Note that linearization with respect to the inputs to the GP will depend on the selected mean $m(\cdot)$ and kernel $k(\cdot,\cdot)$ functions. In general, the above Jacobians can be easily computed via \textit{automatic differentiation} (e.g. Autograd in PyTorch [\cite{paszke2017automatic}]).


\section{Greedy Nonlinear Covariance Steering}
\label{section: greedy cs}
\subsection{Problem Formulation}
Consider the finite-time evolution of the stochastic system \eqref{eq: exact dynamics}. The goal of finite-time covariance steering is to find a control policy that will steer the state of \eqref{eq: exact dynamics} from a given initial distribution with mean $\vmu_0$ and covariance $\Sigma_0$ to a given terminal one with mean and covariance $\vmu_f$ and $\Sigma_f$, respectively, in a finite horizon of $T$ time steps.

A \textit{greedy} approach to finite-horizon covariance steering was presented in \cite{bakolas2020greedy}, where the dynamics \eqref{eq: exact dynamics} are linearized at each time step around the current mean, the \textit{linear} covariance steering problem from the current to the target mean and covariance is solved, and only the first control law is applied -- a model-predictive control approach with a shrinking horizon. However, the exact dynamics in \eqref{eq: exact dynamics} are unknown and cannot be used in the model-based covariance steering algorithm of \cite{bakolas2020greedy}. Instead, the greedy algorithm will be adapted to be used with the approximate, non-parametric GP model that we learned in Section \ref{section: sysid}.

In particular, consider the learned model \eqref{eq: svgp model} for $t = 0,\dots,T-1$, with an initial state $\vz_0$ drawn from a distribution with $\Eb[\vz_0] = \vmu_0$ and $\mathrm{Cov}[\vz_0] = \Sigma_0$, where $\vmu_0 \in \Rb^{n_z}$ and $\Sigma_0 \in \Sb^{++}_{n_z}$ are given. The process noise, $\vw_t$, is assumed to be a sequence of i.i.d. random variables drawn from \eqref{eq: w definition}. Furthermore, $\vz_0$ is conditionally independent of $\vw_t$, for all $t = 0,\dots,T-1$.

Because the identified system in \eqref{eq: svgp model} is nonlinear, there is no guarantee that an initial state drawn from a normal distribution will lead to future states being Gaussian. Therefore, as explained in \cite{bakolas2020greedy}, it is more prudent to talk about steering the nonlinear system mean and covariance \textit{close} to desired quantities rather than steering the state distribution to a goal distribution. 

If we take the class of admissible control policies to be the set of sequences of control laws that are measurable functions of the realization of the current state, the nonlinear covariance steering problem can be formulated as follows:

\begin{problem}[nonlinear covariance steering problem]\label{problem1} 
    \hspace{0.1cm}\\Let $\vmu_0, \vmu_\f \in \Rb^{n_z}$ and $\Sigma_0, \Sigma_\f \in \mathbb{S}^{++}_{n_z}$ be given. Find a control policy $\pi : = \{ \kappa_t(\cdot)\}_{t=0}^{T-1}$ that will steer the system \eqref{eq: svgp model} and, consequently, \eqref{eq: exact dynamics}, from the initial state $\vz_0$ with $\mathbb{E}[\vz_0]=\vmu_0$ and $\mathrm{Cov}[\vz_0] = \Sigma_0$ to a terminal state $\vz_T$ with
    \begin{equation}
        \vmu_T =\vmu_\f,~~~~ (\Sigma_\f - \Sigma_T) \in
        \mathbb{S}^{+}_{n_z}.
    \label{eq: problem_1_constraints}
    \end{equation}
\end{problem}
\textit{Remark:} Given that the system in \eqref{eq: svgp model} is nonlinear, enforcing the equality constraint $\Sigma_T = \Sigma_\f$ would be a difficult task in practice. Following \cite{bakolas2020greedy}, we consider instead the relaxed constraint given in \eqref{eq: problem_1_constraints} according to which, it suffices to achieve a terminal state covariance $ \Sigma_T$ that is ``smaller'' (in the Loewner sense) than $\Sigma_\f$, which corresponds to a situation in which the (desired) terminal mean $\vmu_\f$ will be reached by representative samples of system's trajectories with less uncertainty than the uncertainty corresponding to $\Sigma_\f$.

\subsection{Finite-Horizon Linearized Covariance Steering Problem}
Next, we formulate a linearized covariance steering problem for
the system described by the linearization 
\begin{align}
    \vz_{j+1|t} \approx A_t \vz_{j|t} + B_t \vu_{j|t} + \vd_t,
    \label{eq: linearized dynamics}
\end{align}
of \eqref{eq: svgp model} around the mean state $\vmu_t$ and corresponding (previous) control policy,
\begin{align}
    \{A_t, B_t, \vd_t\} = \mathrm{LIN}_{G}\{\vmu_t, \phi^*_{t|t-1}\left( \{\vmu_i\}_{i=t-1}^t \right)\},
\end{align}
for $j = t,\dots,T-1$. For the latter problem, consider the class $\calU$ of admissible control policies that consist of the sequence of control laws $\{\phi_{j|t}(
\cdot )\}_{j=t}^{T-1}$ that are affine functions of the histories of states, that is,
\begin{equation}
    \phi_{j|t} ( \{\vz_{i|t}\}_{i=t}^j ) = \vupsilon_{j|t} + \sum_{i=t}^{j} K_{j,i|t} \vz_{i|t},
    \label{eq: phi}
\end{equation}
for $j = t,\dots,T-1$. The linearized covariance steering problem at time step $t$ is formulated as follows:

\begin{problem}[$t$-th linearized covariance steering problem] \label{problem2}
    Let $\vmu_t, \vmu_\f\in\mathbb{R}^{n_z}$ and $\Sigma_t, \Sigma_\f \in
    \mathbb{S}^{++}_{n_z}$ be given. Among all admissible control policies
    $\varpi_t : = \{ \phi_{j|t}(\cdot) \}_{j=t}^{T-1} \in
    \calU$, with $\phi_{j|t}(\cdot)$ of the form \eqref{eq: phi}, find a control policy $\varpi_t^*$ that minimizes
    the performance index
    \begin{align}\label{eq:cost}
        J_t(\varpi_t) & :=
        \mathbb{E}\big[\sum_{j = t}^{T-1} \phi_{j|t}( \{\vz_{i|t}\}_{i=t}^j )^\top \phi_{j|t}( \{\vz_{i|t}\}_{i=t}^j )
        \big]
    \end{align}
    subject to the recursive dynamic constraints \eqref{eq: linearized dynamics} and the boundary conditions
    \begin{subequations}
        \begin{align}
            \Eb[\vz_{t|t}] & = \vmu_t, & \mathrm{Cov}[\vz_{t|t}] = \Sigma_t,\\
            \Eb[\vz_{T|t}] & =\vmu_\f, & (\Sigma_\f - \mathrm{Cov}[\vz_{T|t}] ) \in
            \mathbb{S}^{+}_{n_z}.
        \end{align}
    \end{subequations}
\end{problem}

The choice of the performance index ensures that the control input will have finite energy, without excessive actuation. Note that the terminal positive semi-definite constraint $(\Sigma_\f - \mathrm{Cov}[\vz_{T|t}]) \in \mathbb{S}^{+}_{n_z}$ differentiates Problem \ref{problem2} from the standard linear quadratic Gaussian (LQG) problem. Although no state or input constraints are considered in this formulation, the optimization-based solution presented here and in \cite{bakolas2020greedy} is applicable to the general problem formulation that includes such constraints (refer to \cite{p:BAKOLAS2018} for more details).

Note that Problem \ref{problem2} can be formulated as a convex \textit{semi-definite program} (SDP) and, thus, can be solved efficiently using any available conic solver. The formulation of the SDP is ommitted in the interest of space, but can be found in \cite{p:BAKOLAS2018} and \cite{bakolas2020greedy}.

For compactness, denote the solution to the $t$-th linearized covariance steering problem as
\begin{align}
    \{ \phi^*_{j|t}(\cdot) \}_{j=t}^{T-1} = \mathrm{LCS}_{t,T}\{A_t, B_t, \vd_t, \vmu_t, \Sigma_t, \vmu_\f, \Sigma_\f\}.
\end{align}

\subsection{Gaussian Process-Based Unscented Transform for Uncertainty Propagation}
\label{subsection: unscented transform}
Let $\pi = \{\kappa_t(\cdot)\}_{t=0}^{T-1}$ be an admissible control policy for Problem \ref{problem1}. Then, the closed-loop dynamics become
\begin{equation}
    \vz_{t+1} = G(\vz_t, \kappa_t(\vz_t)) + \vw_t.
    \label{eq: closed loop gp}
\end{equation}

The mean and covariance of the uncertain state of the nonlinear system described by \eqref{eq: closed loop gp} is propagated using the \textit{unscented transform} [\cite{p:julier2002,p:julier2004}]. To this aim, assume that the mean $\vmu_t:= \mathbb{E}[\vz_t]$ and covariance $\Sigma_t:=\mathrm{Cov}[\vz_t]$ of the state of \eqref{eq: svgp model} (or estimates of these quantities) are known at time step $t$. 

First, we compute $2n_z+1$ deterministic points, $\sigma_t^{(i)},\ i = 1,\dots 2n_z+1$, which are also known as \textit{sigma points}, according to \cite{p:julier2002, p:julier2004}. Then, to each sigma point, we associate a pair of gains $(\gamma_t^{(i)}, \delta_t^{(i)})$, according to \cite{p:julier2004, p:wan2000}. Subsequently, the sigma points $\{\sigma_t^{(i)}\}_{i=1}^{2n_z+1}$ are propagated to the next time step to obtain a new set of points $\{ \hat{\sigma}_{t+1}^{(i)} \}_{i=1}^{2n_z+1}$, where
\begin{equation}\label{eq:sigma2}
    \hat{\sigma}_{t+1}^{(i)} = G(\sigma_t^{(i)}, \kappa_t(\sigma_t^{(i)})),~~~i = 1,\dots,2n_z.
\end{equation}
Using this new point-set, one can approximate the (predicted) state mean and covariance at time step $t+1$ as
\begin{subequations}
    \begin{align}\label{eq:Xiupdate}
        \hat{\vmu}_{t+1} & = \sum_{i=0}^{2n_z} \gamma_t^{(i)} \hat{\sigma}_{t+1}^{(i)},\\
        \hat{\Sigma}_{t+1} & = \sum_{i=0}^{2n_z} \delta_t^{(i)} (
        \hat{\sigma}_{t+1}^{(i)} - \hat{\vmu}_{t+1} ) ( \hat{\sigma}_{t+1}^{(i)} -
        \hat{\vmu}_{t+1} )^\top + W_t. \label{eq:Xiupdate2}
    \end{align}
\end{subequations}
Similar to \cite{ko2007gp}, we set $W_t = \mathrm{Cov}[\vw_t] = \Sigma_f([\vz_t; \kappa_t(\vz_t)], [\vz_t; \kappa_t(\vz_t)]) + \sigma_{\epsilon}^2 \in \mathbb{S}^{+}_{n_z}$ as the process noise covariance. Notice that $W_t$ captures both the noise in the system as well as the model uncertainties resulting from the lack of training data points used in the learning phase.

\subsection{Greedy Nonlinear Covariance Steering for Gaussian Process Predictive Models}
Now we have all the tools necessary to extend the \textit{greedy nonlinear covariance steering} algorithm of \cite{bakolas2020greedy} to Gaussian process predictive models. The greedy algorithm consists of three main steps. Consider the time step $t$, where $t=0,\dots,T-1$, and assume that estimates of the state mean, $\hat{\vmu}_t$, the state covariance, $\hat{\Sigma}_t$, as well as the input mean $\hat{\vnu}_t$, are known (starting from $\hat{\vmu}_0 = \vmu_0$, $\hat{\Sigma}_0 = \Sigma_0$, and $\hat{\vnu}_0 = 0$).

The \textit{first step} is to linearize \eqref{eq: svgp model} around $(\hat{\vmu}_t, \hat{\vnu}_t)$:
\begin{align}
    \{A_t, B_t, \vd_t\} = \mathrm{LIN}_{G}\{\hat{\vmu}_t, \hat{\vnu}_t\},
\end{align}
where $\hat{\vnu}_t = \phi^*_{t|t-1}(\{\hat{\vmu}_i\}_{i=t-1}^t)$. The linearization will have to be updated at each time step $t$ since the estimates $\hat{\vmu}_t$ and $\hat{\vnu}_t$ will also be updated.

The \textit{second step} is to solve the $t$\textit{-th linearized covariance steering} problem and compute the feedback control policy that solves Problem \ref{problem2}. The latter problem is solved using the linearized model $\{A_t,B_t,\vd_t\}$ obtained in the first step and the estimates of the predicted mean and covariance $(\hat{\vmu}_t, \hat{\Sigma}_t)$ at time step $t$:
\begin{align}
    \{ \phi^*_{j|t}(\cdot) \}_{j=t}^{T-1} = \mathrm{LCS}_{t,T}\{A_t, B_t, \vd_t, \hat{\vmu}_t, \hat{\Sigma}_t, \vmu_\f, \Sigma_\f\}.
\end{align}
The computation of $\{ \phi^*_{j|t}(\cdot) \}_{j=t}^{T-1}$ can be done in real-time by means of robust and efficient convex optimization techniques [\cite{p:BAKOLAS2018}]. 
 
From $\{ \phi^*_{j|t}(\cdot) \}_{j=t}^{T-1}$, we extract only the first control law,
\begin{align*}
    \kappa_t(\vz_t) := \phi^*_{t|t}(\vz_t) = \vupsilon^*_{t|t} +
    K^*_{t|t} \vz_t,
\end{align*}
where $\vz$ is the state of the original nonlinear system.
The one-time-step transition map for the closed-loop dynamics based on information available at time step $t$ is then described by
\begin{align}
    \vz_{t+1} &= G(\vz_t, \kappa_t(\vz_t)) + \vw_t.
\end{align}

In the \textit{third step}, the estimates $(\hat{\vmu}_{t+1}, \hat{\Sigma}_{t+1})$ are computed by propagating the mean $\hat{\vmu}_t$ and covariance $\hat{\Sigma}_t$ of the closed-loop system to the next time step. The new mean and covariance, i.e., $\hat{\vmu}_{t+1}$ and $\hat{\Sigma}_{t+1}$, are computed using the GP-based unscented transform described in Section \ref{subsection: unscented transform}. We write
\begin{equation}
    (\hat{\vmu}_{t+1}, \hat{\Sigma}_{t+1}) :=
    \mathrm{UT}\{ \hat{\vmu}_t,\hat{\Sigma}_t;
    G(\cdot, \cdot), \kappa_t(\cdot) \}.
\end{equation}

The three steps of the greedy covariance steering algorithm are repeated for all time steps $t = 0,\dots,T-1$. At the end of the process, the predicted approximations of the state mean and covariance should be sufficiently close to their corresponding goal quantities. The output of this iterative process is a control policy $\pi :=\{\kappa_t(\vz) \}_{t=0}^{T-1}$ that solves Problem \ref{problem1}.


\section{Numerical Results}
\label{section: results}
In this section, the basic ideas of this paper are illustrated in numerical simulations. In particular, consider the following stochastic nonlinear system:
\begin{subequations}\label{eq: car dynamics}
    \begin{align}
        s_{x,t+1} & = s_{x,t} + v_t \tau \cos{\theta_t} + \epsilon^{s_x}_t,\\
        s_{y,t+1} & = s_{y,t} + v_t \tau \sin{\theta_t} + \epsilon^{s_y}_t,\\
        \theta_{t+1} & = \theta_t + u_t^{\theta} v_t \tau + \epsilon^{\theta}_t,\\
        v_{t+1} & = v_t + u_t^{v} \tau  + \epsilon^{v}_t
    \end{align}
\end{subequations}
which is a discrete-time realization of a unicycle car model with state $\vz = \begin{bmatrix} s_x & s_y & \theta & v\end{bmatrix}^\top$ and input $\vu = \begin{bmatrix} u^{\theta} & u^v\end{bmatrix}^\top$. In order to train the SVGP model, we assume that a ``black-box'' simulator of the dynamics \eqref{eq: car dynamics} is available for sampling full-state transitions $\vy_i$ for given states $\vz_i$ and inputs $\vu_i$ (see \eqref{eq: observations} and \eqref{eq: inputs}, respectively). 
We run the simulator and collect a set of training data points that are then used to learn a stochastic variational GP model of the system dynamics, as presented in Section \ref{section: sysid}. 
Then, the learned model is used to steer the mean and covariance of the state of the underlying system from a given initial distribution to a prescribed terminal one. For our simulations, we consider \eqref{eq: car dynamics} with time step $\tau = 0.05$, while the white noise $\epsilon_t \sim \calN\left(\epsilon_t \mid 0, \mathrm{diag}(\begin{bmatrix} 0.02 & 0.02 & 0.04 & 0.04 \end{bmatrix}^2\right)$.

\subsubsection{Training: } The SVGP model with zero mean, $m(\vx) = 0$, squared-exponential kernel, $$k(\vx, \vx') = \sigma_f^2 \exp{\left(-\frac{1}{2}\left(\vx - \vx'\right)^\top L^{-1}\left(\vx - \vx'\right)\right)}$$ with separate length scales $L = \mathrm{diag}(l_1^2,\dots,l_n^2)$ for each input dimension, and $M = 256$ inducing locations is setup using GPyTorch [\cite{gardner2018gpytorch}]. A set of $N=9000$ data points $(\vx_i, \vy_i)$ are collected from randomly sampled states $\vz$ and inputs $\vu$ between $\vz_{min} = \begin{bmatrix} -20 & -20 & -6\pi & -10\end{bmatrix}^\top$, $\vu_{min} = \begin{bmatrix} -20 & -20\end{bmatrix}^\top$ and  $\vz_{max} = \begin{bmatrix} 20 & 20 & 6\pi & 20\end{bmatrix}^\top$, $\vu_{max} = \begin{bmatrix} 20 & 20\end{bmatrix}^\top$, respectively. The variational parameters $\vZ$, $\vm$, and $\vS$, along with the hyperparameters $\mathbf{\Theta} = \{\theta_m, \theta_t, \sigma_{\epsilon} \}$, are optimized by minimizing the negative lower bound, $-\calL$, using the Adam optimizer [\cite{kingma2014adam}]. 

\subsubsection{GP-Based Greedy Covariance Steering: } Assume an initial state $\vz_0$ drawn from a distribution with mean $\vmu_0 = \begin{bmatrix} 0 & 0 & 0 & 1\end{bmatrix}^\top$ and covariance $\Sigma_0 = \mathrm{diag}(\begin{bmatrix} [0.1 & 0.2 & 0.1 & 0.1] \end{bmatrix}^2$. The \textit{target}  terminal state mean and covariance are taken to be $\vmu_\f = \begin{bmatrix} 1 & 2 & 0 & 1\end{bmatrix}^\top$ and $\Sigma_\f = \mathrm{diag}(\begin{bmatrix} [0.1 & 0.05 & 0.05 & 0.05] \end{bmatrix}^2$, respectively. The greedy covariance steering algorithm of Section \ref{section: greedy cs} is run with the identified SVGP model with a time horizon of $T=30$ time steps. The Jacobians \eqref{eq: jacobians} for the model linearizations are computed using automatic differentiation in PyTorch [\cite{paszke2017automatic}]. For comparison, the \textit{exact} model \eqref{eq: car dynamics} is used with the greedy algorithm of \cite{bakolas2020greedy} for the same initial and target distributions.

\begin{figure}
     \centering
     \begin{subfigure}[t]{0.49\linewidth}
         \centering
         \includegraphics[scale=0.36,trim={3.5cm 0 5.5cm 0},clip]{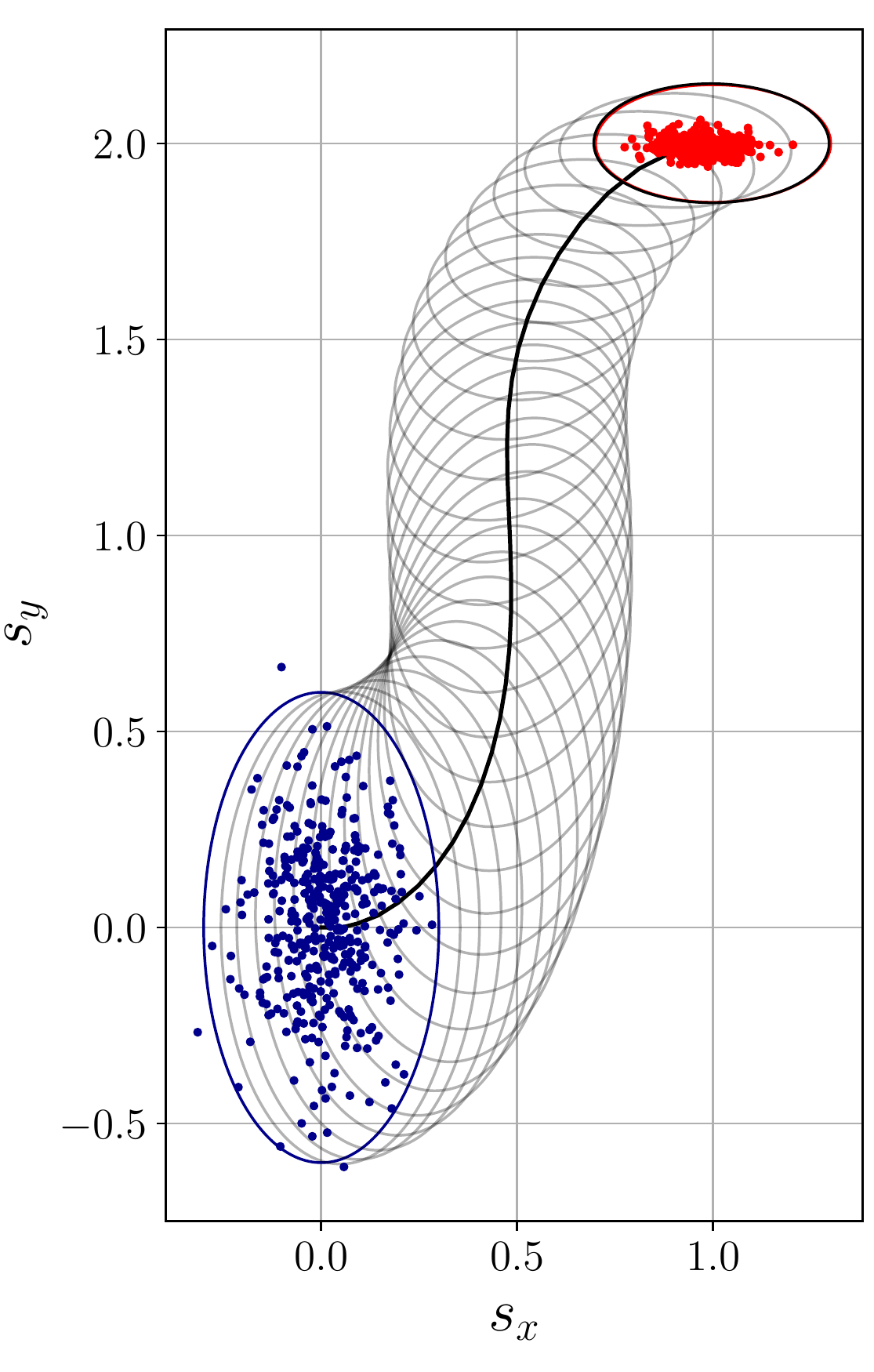}
         \caption{SVGP Model}
         \label{fig: position gp}
     \end{subfigure}%
     \begin{subfigure}[t]{0.49\linewidth}
         \centering
         \includegraphics[scale=0.36,trim={3.5cm 0 5.5cm 0},clip]{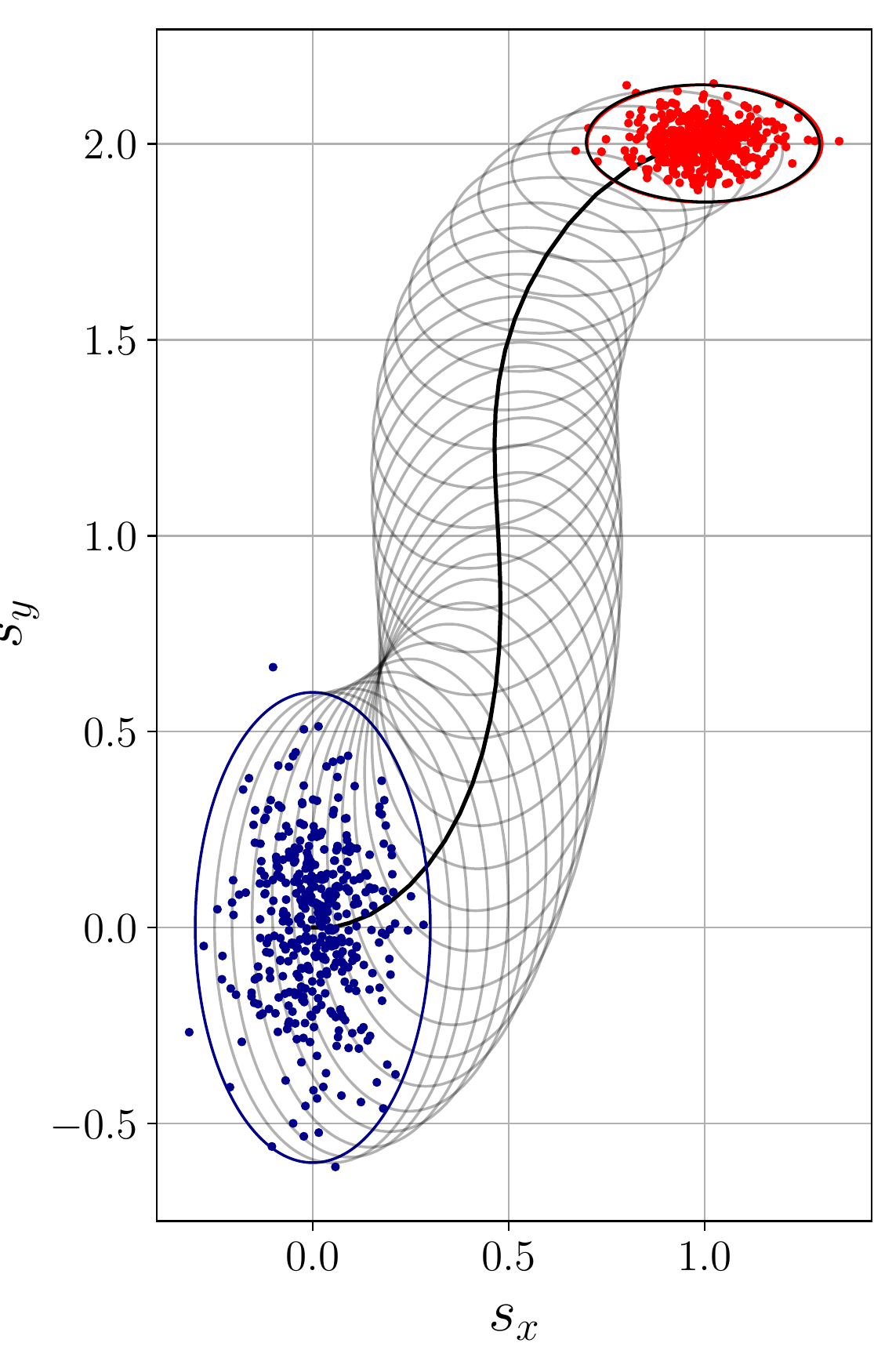}
         \caption{Exact Model}
         \label{fig: position exact}
     \end{subfigure}
        \caption{Position ($s_x$, $s_y$) uncertainties. Blue: initial distribution. Red: target distribution. Black: actual terminal distribution. Gray: intermediate distributions. Solid black line: mean trajectory.}
        \label{fig: positions}
\end{figure}

With the chosen target covariance, the goal is to shrink the uncertainty in the coordinate $s_y$, the angle $\theta$, and the velocity $v$, while retaining the uncertainty in $s_x$. The mean trajectory of $(s_x,s_y)$ along with the corresponding $3-\sigma$ uncertainty ellipsoids are plotted in Fig. \ref{fig: positions}. The GP-based algorithm (Fig. \ref{fig: position gp}) is compared with the greedy algorithm that uses the exact (analytic) model of \eqref{eq: car dynamics} (Fig. \ref{fig: position exact}). The first thing we notice is that the uncertainty predicted by the SVGP model is overestimated. This is expected, since GPs provide conservative uncertainty estimates that include not only the uncertainties due to process noise, but also due to modeling errors. The terminal distribution estimated by the Unscented Transform (solid black ellipsoid) reaches almost perfectly the desired one (red ellipsoid). However, since the uncertainties provided by the SVGP model are conservative, the actual distribution -- which is visualized as particles (red) from $400$ Monte Carlo realizations -- is shrinked compared to the target one. In comparison, the estimated terminal distribution for the exact model reaches almost perfectly the target one (see red particles in Fig. \ref{fig: position exact}). Thus, the use of SVGP model results in a more \textit{cautious} covariance steering. 

\begin{figure}
     \centering
     \begin{subfigure}[t]{0.49\linewidth}
         \centering
         \includegraphics[scale=0.4,trim={0 0 0 0},clip]{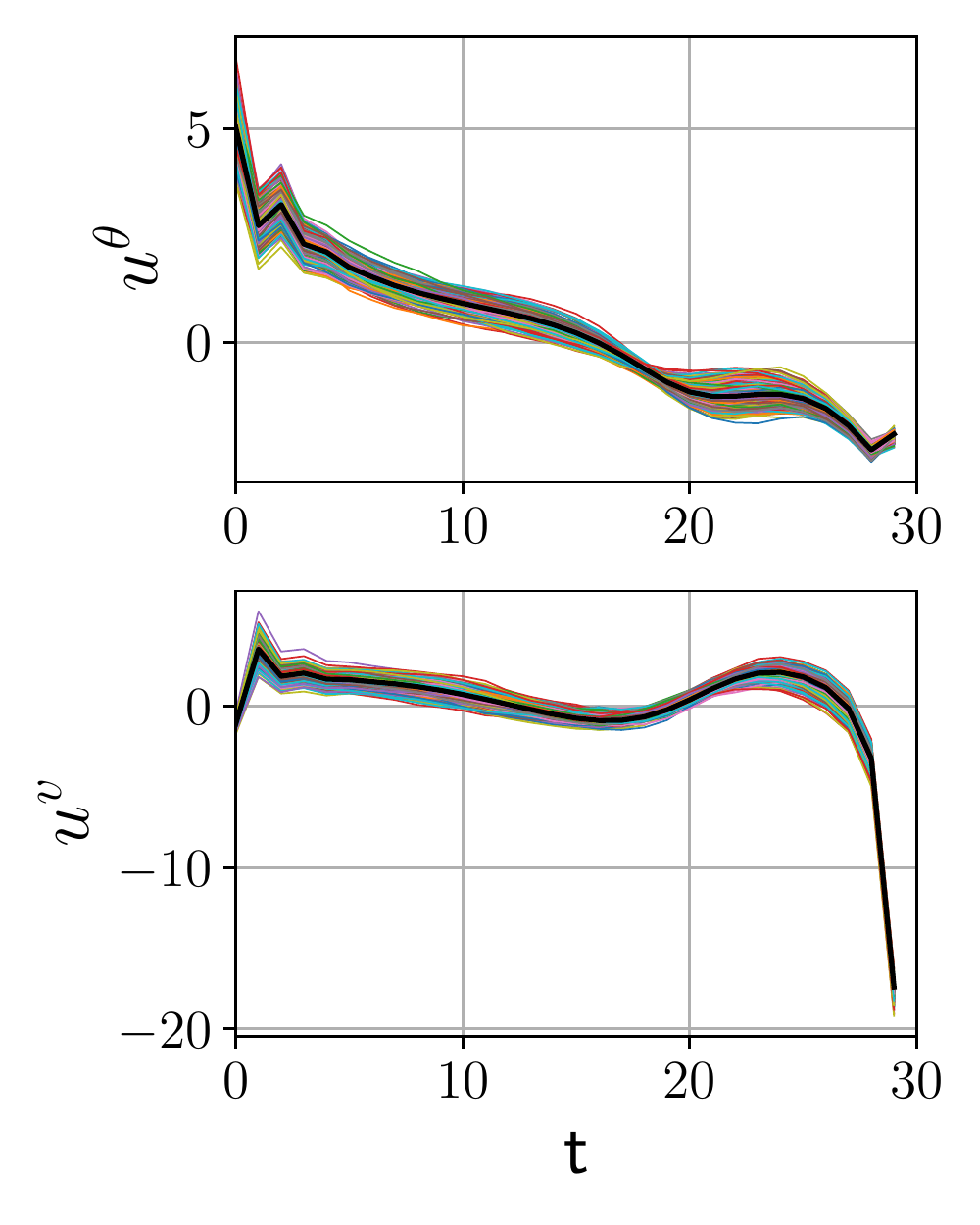}
         \caption{SVGP Model}
         \label{fig: inputs gp}
     \end{subfigure}%
     \begin{subfigure}[t]{0.49\linewidth}
         \centering
         \includegraphics[scale=0.4,trim={0 0 0 0},clip]{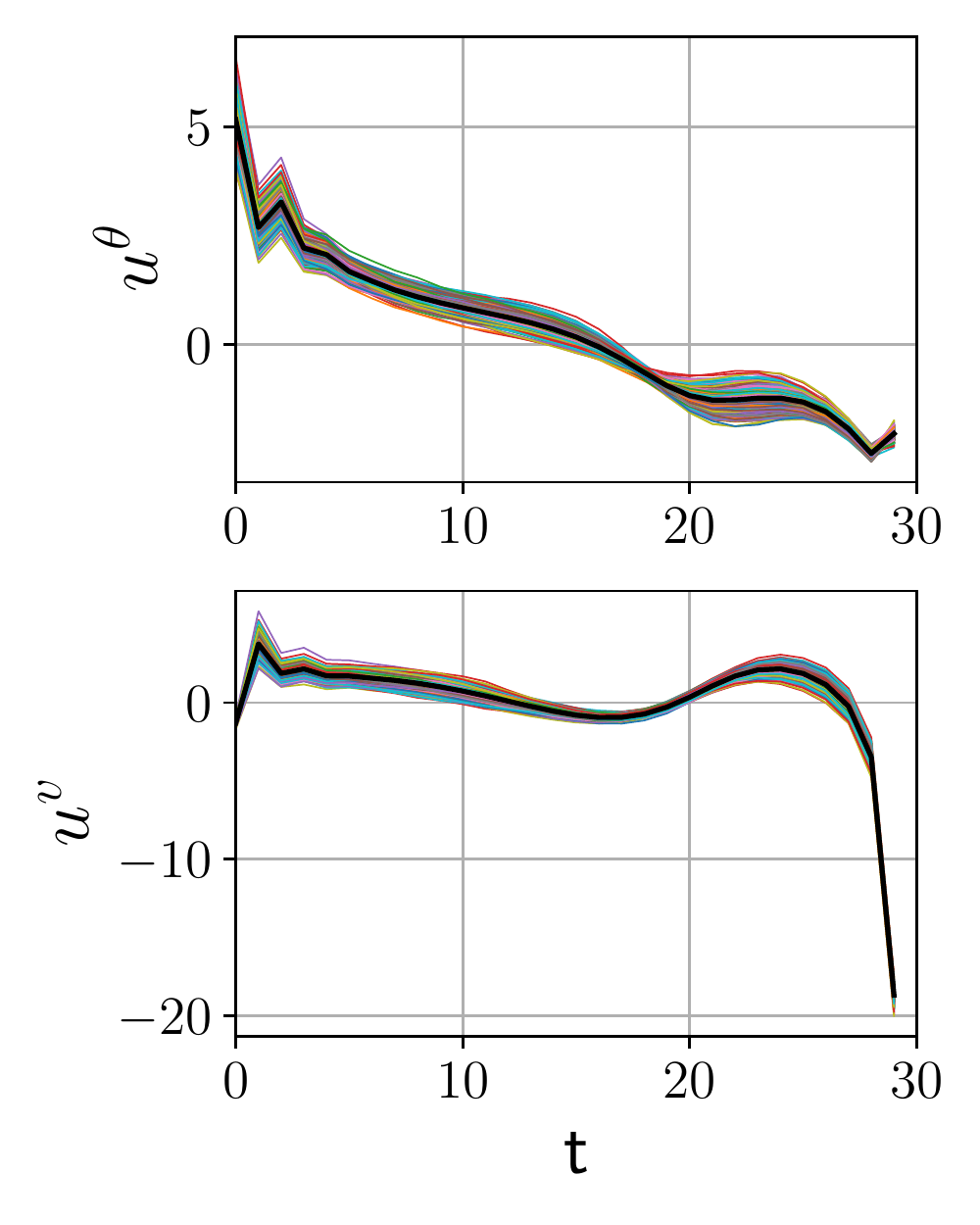}
         \caption{Exact Model}
         \label{fig: inputs exact}
     \end{subfigure}
        \caption{Inputs $\vu_t$. Solid black: mean input $\vnu_t$. Colored: inputs of Monte Carlo realizations.}
        \label{fig: inputs}
\end{figure}

\begin{figure}
     \centering
     \begin{subfigure}[t]{0.49\linewidth}
         \centering
         \includegraphics[scale=0.4,trim={0 0 0 0},clip]{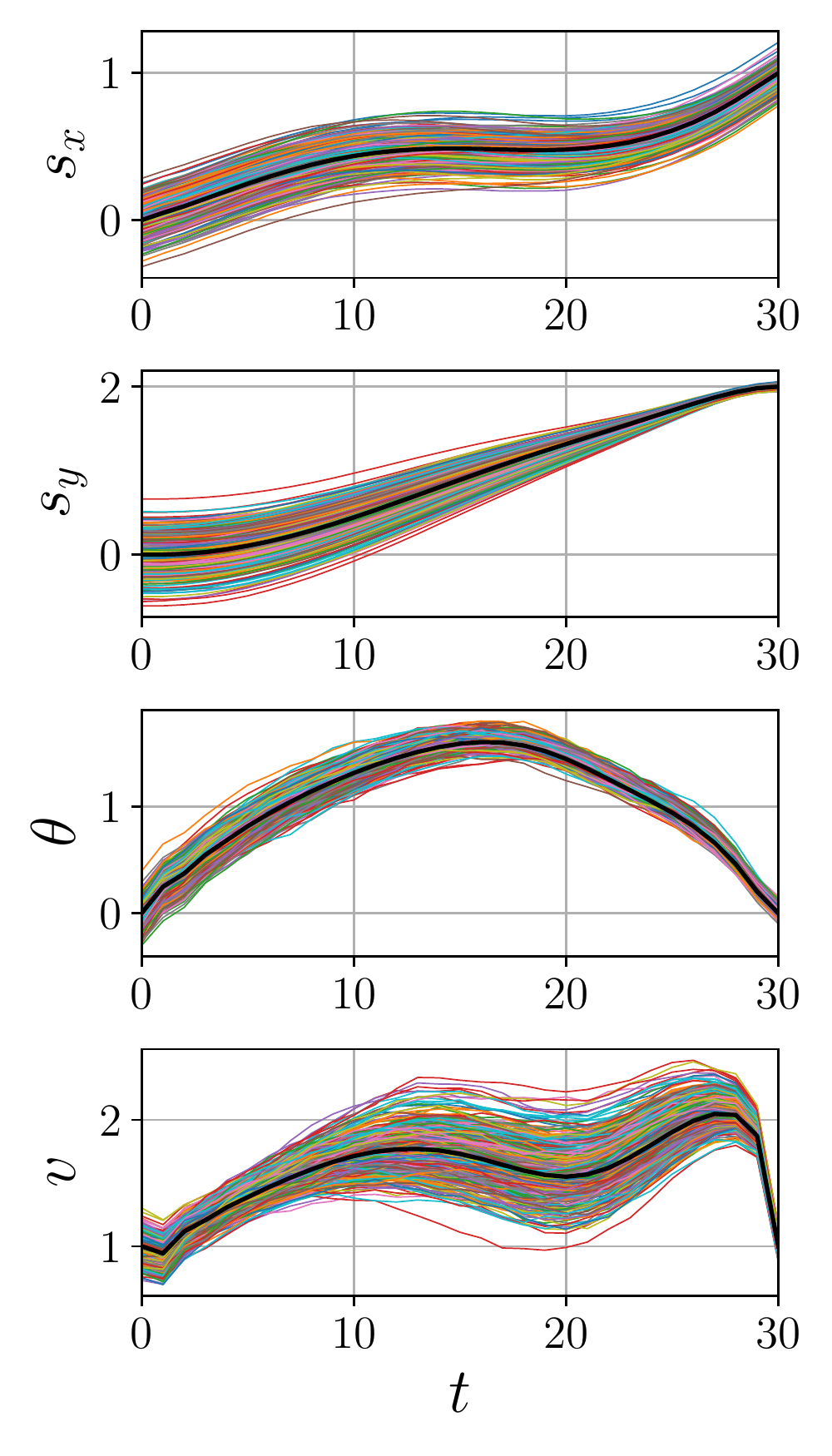}
         \caption{SVGP Model}
         \label{fig: states gp}
     \end{subfigure}%
     \begin{subfigure}[t]{0.49\linewidth}
         \centering
         \includegraphics[scale=0.4,trim={0 0 0 0},clip]{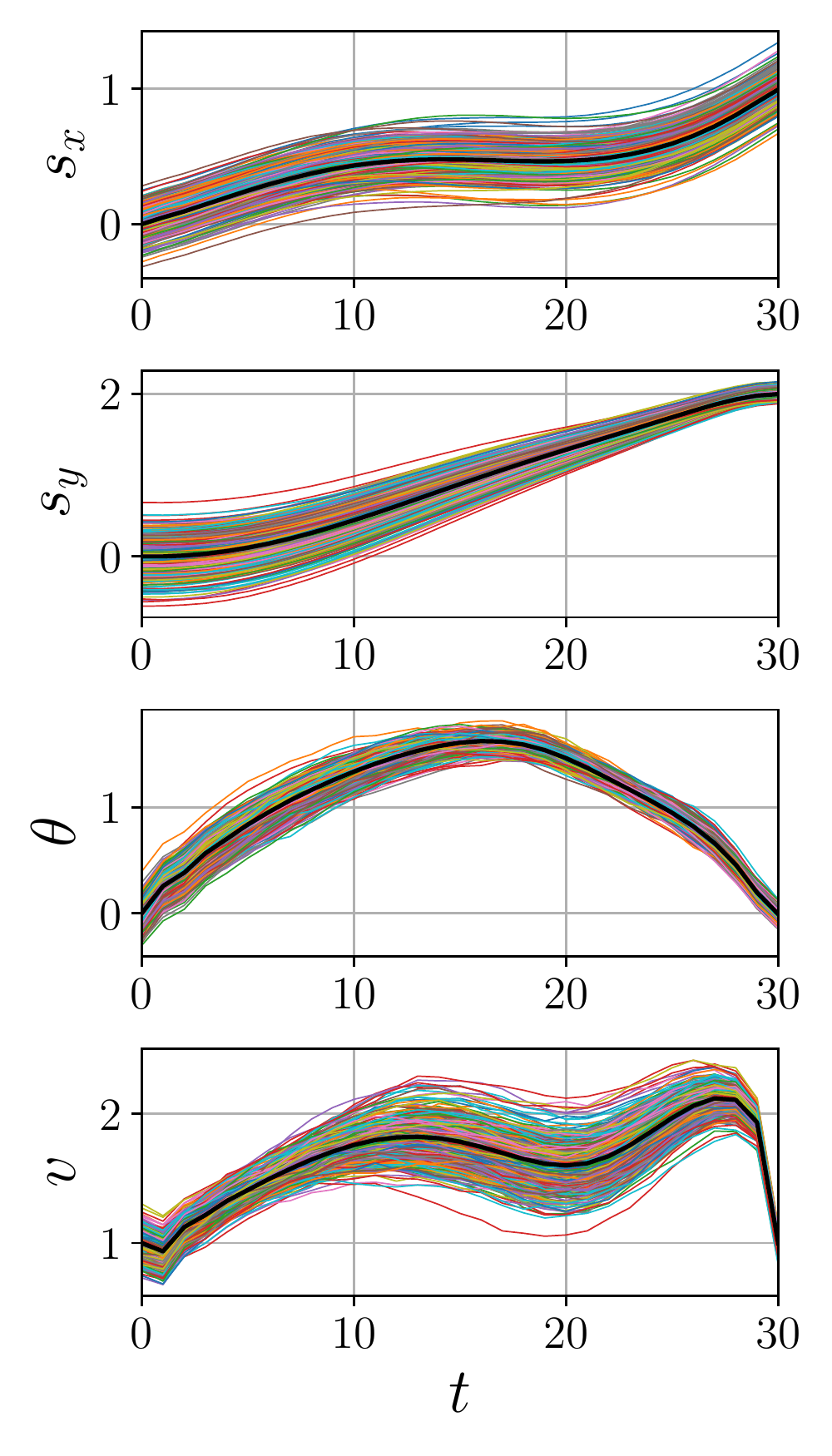}
         \caption{Exact Model}
         \label{fig: states exact}
     \end{subfigure}
        \caption{States $\vz_t$. Solid black: mean state $\vz_t$. Colored: states of Monte Carlo realizations.}
        \label{fig: states}
\end{figure}

\section{Conclusion}
\label{section: conclusions}
In this work, a greedy covariance steering algorithm that uses scalable Gaussian process predictive models for discrete-time stochastic nonlinear systems with unknown dynamics has been proposed. First, a non-parametric predictive model is learned from a set of training data points using stochastic variational Gaussian process regression. Then, a set of linearized covariance steering problems is solved and the mean and covariance of the closed-loop system is predicted using the unscented transform. This work has considered the case of perfect full-state information. However, more practical cases, such as that of incomplete information, where the states of the system have to be estimated from partial measurements, will be explored in the future by the authors.

\bibliography{biblio}             

\begin{thebibliography}{33}
\providecommand{\natexlab}[1]{#1}
\providecommand{\url}[1]{\texttt{#1}}
\providecommand{\urlprefix}{URL }
\expandafter\ifx\csname urlstyle\endcsname\relax
  \providecommand{\doi}[1]{doi:\discretionary{}{}{}#1}\else
  \providecommand{\doi}{doi:\discretionary{}{}{}\begingroup
  \urlstyle{rm}\Url}\fi

\bibitem[{Bakolas(2016)}]{p:bakcdc16}
Bakolas, E. (2016).
\newblock Optimal covariance control for discrete-time stochastic linear
  systems subject to constraints.
\newblock In \emph{IEEE CDC (2016)}, 1153--1158.

\bibitem[{Bakolas(2017)}]{p:EBACC17}
Bakolas, E. (2017).
\newblock Covariance control for discrete-time stochastic linear systems with
  incomplete state information.
\newblock In \emph{2017 American Control Conference (ACC)}, 432--437.

\bibitem[{Bakolas(2018)}]{p:BAKOLAS2018}
Bakolas, E. (2018).
\newblock Finite-horizon covariance control for discrete-time stochastic linear
  systems subject to input constraints.
\newblock \emph{Automatica}, 91, 61--68.

\bibitem[{{Bakolas}(2019)}]{p:bakTAC2019}
{Bakolas}, E. (2019).
\newblock Dynamic output feedback control of the liouville equation for
  discrete-time siso linear systems.
\newblock \emph{IEEE Transactions on Automatic Control}, 64(10), 4268--4275.

\bibitem[{Bakolas and Tsolovikos(2020)}]{bakolas2020greedy}
Bakolas, E. and Tsolovikos, A. (2020).
\newblock Greedy finite-horizon covariance steering for discrete-time
  stochastic nonlinear systems based on the unscented transform.
\newblock In \emph{ACC (2020)}, 3595--3600.

\bibitem[{Caluya and Halder(2019)}]{p:caluya2019}
Caluya, K.F. and Halder, A. (2019).
\newblock Finite horizon density control for static state feedback linearizable
  systems.
\newblock \emph{arXiv preprint arXiv:1904.02272}.

\bibitem[{Chen et~al.(2016{\natexlab{a}})Chen, Georgiou, and
  Pavon}]{p:georgiou15A}
Chen, Y., Georgiou, T., and Pavon, M. (2016{\natexlab{a}}).
\newblock Optimal steering of a linear stochastic system to a final probability
  distribution, \textsc{P}art \textsc{I}.
\newblock \emph{IEEE Trans. on Autom. Control}, 61(5), 1158 -- 1169.

\bibitem[{Chen et~al.(2016{\natexlab{b}})Chen, Georgiou, and
  Pavon}]{p:georgiou15B}
Chen, Y., Georgiou, T., and Pavon, M. (2016{\natexlab{b}}).
\newblock Optimal steering of a linear stochastic system to a final probability
  distribution, \textsc{P}art \textsc{II}.
\newblock \emph{IEEE Trans. Autom. Control}, 61(5), 1170--1180.

\bibitem[{Gardner et~al.(2018)Gardner, Pleiss, Bindel, Weinberger, and
  Wilson}]{gardner2018gpytorch}
Gardner, J.R., Pleiss, G., Bindel, D., Weinberger, K.Q., and Wilson, A.G.
  (2018).
\newblock {GPyTorch}: Blackbox matrix-matrix gaussian process inference with
  {GPU} acceleration.
\newblock \emph{arXiv preprint arXiv:1809.11165}.

\bibitem[{Goldshtein and Tsiotras(2017)}]{p:PT2017}
Goldshtein, M. and Tsiotras, P. (2017).
\newblock Finite-horizon covariance control of linear time-varying systems.
\newblock In \emph{IEEE CDC (2017)}, 3606--3611.

\bibitem[{Grigoriadis and Skelton(1997)}]{p:Grig97}
Grigoriadis, K.M. and Skelton, R.E. (1997).
\newblock Minimum-energy covariance controllers.
\newblock \emph{Automatica}, 33(4), 569--578.

\bibitem[{Grimes et~al.(2006)Grimes, Chalodhorn, and Rao}]{grimes2006dynamic}
Grimes, D.B., Chalodhorn, R., and Rao, R.P. (2006).
\newblock Dynamic imitation in a humanoid robot through nonparametric
  probabilistic inference.
\newblock In \emph{Robotics: science and systems}, 199--206. Cambridge, MA.

\bibitem[{Hensman et~al.(2013)Hensman, Fusi, and
  Lawrence}]{hensman2013gaussian}
Hensman, J., Fusi, N., and Lawrence, N.D. (2013).
\newblock Gaussian processes for big data.
\newblock \emph{arXiv preprint arXiv:1309.6835}.

\bibitem[{Hoffman et~al.(2013)Hoffman, Blei, Wang, and
  Paisley}]{hoffman2013stochastic}
Hoffman, M.D., Blei, D.M., Wang, C., and Paisley, J. (2013).
\newblock Stochastic variational inference.
\newblock \emph{Journal of Machine Learning Research}, 14(5).

\bibitem[{Hotz and Skelton(1987)}]{p:skeltonIJC}
Hotz, A. and Skelton, R.E. (1987).
\newblock Covariance control theory.
\newblock \emph{Int. J. Control}, 16(1), 13--32.

\bibitem[{{Julier} and {Uhlmann}(2004)}]{p:julier2004}
{Julier}, S.J. and {Uhlmann}, J.K. (2004).
\newblock Unscented filtering and nonlinear estimation.
\newblock \emph{Proc. IEEE}, 92(3), 401--422.

\bibitem[{Julier(2002)}]{p:julier2002}
Julier, S.J. (2002).
\newblock The scaled unscented transformation.
\newblock In \emph{ACC (2002)}, volume~6, 4555--4559. IEEE.

\bibitem[{Kingma and Ba(2014)}]{kingma2014adam}
Kingma, D.P. and Ba, J. (2014).
\newblock Adam: A method for stochastic optimization.
\newblock \emph{arXiv preprint arXiv:1412.6980}.

\bibitem[{Ko and Fox(2009)}]{ko2009gp}
Ko, J. and Fox, D. (2009).
\newblock \textsc{GP-B}ayesfilters: Bayesian filtering using gaussian process
  prediction and observation models.
\newblock \emph{Autonomous Robots}, 27(1), 75--90.

\bibitem[{Ko et~al.(2007{\natexlab{a}})Ko, Klein, Fox, and
  Haehnel}]{ko2007gaussian}
Ko, J., Klein, D.J., Fox, D., and Haehnel, D. (2007{\natexlab{a}}).
\newblock Gaussian processes and reinforcement learning for identification and
  control of an autonomous blimp.
\newblock In \emph{ICRA (2007)}, 742--747.

\bibitem[{Ko et~al.(2007{\natexlab{b}})Ko, Klein, Fox, and Haehnel}]{ko2007gp}
Ko, J., Klein, D.J., Fox, D., and Haehnel, D. (2007{\natexlab{b}}).
\newblock \textsc{GP-UKF}: Unscented kalman filters with gaussian process
  prediction and observation models.
\newblock In \emph{IROS (2007)}, 1901--1907.

\bibitem[{{Mesbah} et~al.(2014){Mesbah}, {Streif}, {Findeisen}, and
  {Braatz}}]{p:mesbah2014}
{Mesbah}, A., {Streif}, S., {Findeisen}, R., and {Braatz}, R.D. (2014).
\newblock Stochastic nonlinear model predictive control with probabilistic
  constraints.
\newblock In \emph{ACC (2014)}, 2413--2419.

\bibitem[{Mukadam et~al.(2016)Mukadam, Yan, and Boots}]{mukadam2016gaussian}
Mukadam, M., Yan, X., and Boots, B. (2016).
\newblock Gaussian process motion planning.
\newblock In \emph{ICRA (2016)}, 9--15.

\bibitem[{Pan and Theodorou(2014)}]{pan2014probabilistic}
Pan, Y. and Theodorou, E. (2014).
\newblock Probabilistic differential dynamic programming.
\newblock In \emph{NIPS (2014)}, 1907--1915.

\bibitem[{Pan and Theodorou(2015)}]{pan2015data}
Pan, Y. and Theodorou, E.A. (2015).
\newblock Data-driven differential dynamic programming using gaussian
  processes.
\newblock In \emph{ACC (2015)}, 4467--4472.

\bibitem[{Paszke et~al.(2017)Paszke, Gross, Chintala, Chanan, Yang, DeVito,
  Lin, Desmaison, Antiga, and Lerer}]{paszke2017automatic}
Paszke, A., Gross, S., Chintala, S., Chanan, G., Yang, E., DeVito, Z., Lin, Z.,
  Desmaison, A., Antiga, L., and Lerer, A. (2017).
\newblock Automatic differentiation in {PyTorch}.

\bibitem[{Quinonero-Candela and Rasmussen(2005)}]{quinonero2005unifying}
Quinonero-Candela, J. and Rasmussen, C.E. (2005).
\newblock A unifying view of sparse approximate gaussian process regression.
\newblock \emph{The Journal of Machine Learning Research}, 6, 1939--1959.

\bibitem[{Rasmussen(2003)}]{rasmussen2003gaussian}
Rasmussen, C.E. (2003).
\newblock Gaussian processes in machine learning.
\newblock In \emph{Summer School on Machine Learning}, 63--71. Springer.

\bibitem[{Ridderhof et~al.(2019)Ridderhof, Okamoto, and
  Tsiotras}]{p:ridderhof2019}
Ridderhof, J., Okamoto, K., and Tsiotras, P. (2019).
\newblock Nonlinear uncertainty control with iterative covariance steering.
\newblock In \emph{CDC (2020)}, 3484--3490.

\bibitem[{Sehr and Bitmead(2017)}]{p:SEHR2017}
Sehr, M.A. and Bitmead, R.R. (2017).
\newblock Particle model predictive control: Tractable stochastic nonlinear
  output-feedback \textsc{MPC}.
\newblock In \emph{20th IFAC World Congress}, 15361 -- 15366.

\bibitem[{Titsias(2009)}]{titsias2009variational}
Titsias, M. (2009).
\newblock Variational learning of inducing variables in sparse gaussian
  processes.
\newblock In \emph{Artificial intelligence and statistics}, 567--574. PMLR.

\bibitem[{Wan and Van Der~Merwe(2000)}]{p:wan2000}
Wan, E.A. and Van Der~Merwe, R. (2000).
\newblock The unscented \textsc{K}alman filter for nonlinear estimation.
\newblock In \emph{Proceedings of the IEEE 2000 Adaptive Systems for Signal
  Processing, Communications, and Control Symposium}, 153--158.

\bibitem[{Xu and Skelton(1992)}]{p:skeltonTAC}
Xu, J.H. and Skelton, R.E. (1992).
\newblock An improved covariance assignment theory for discrete systems.
\newblock \emph{IEEE Trans. Autom. Control}, 37(10), 1588--1591.

\end{thebibliography}


\begin{thebibliography}{4}
\providecommand{\natexlab}[1]{#1}
\providecommand{\url}[1]{\texttt{#1}}
\providecommand{\urlprefix}{URL }
\expandafter\ifx\csname urlstyle\endcsname\relax
  \providecommand{\doi}[1]{doi:\discretionary{}{}{}#1}\else
  \providecommand{\doi}{doi:\discretionary{}{}{}\begingroup
  \urlstyle{rm}\Url}\fi

\bibitem[{Able(1956)}]{Abl:56}
Able, B. (1956).
\newblock Nucleic acid content of microscope.
\newblock \emph{Nature}, 135, 7--9.

\bibitem[{Able et~al.(1954)Able, Tagg, and Rush}]{AbTaRu:54}
Able, B., Tagg, R., and Rush, M. (1954).
\newblock Enzyme-catalyzed cellular transanimations.
\newblock In A.~Round (ed.), \emph{Advances in Enzymology}, volume~2, 125--247.
  Academic Press, New York, 3rd edition.

\bibitem[{Keohane(1958)}]{Keo:58}
Keohane, R. (1958).
\newblock \emph{Power and Interdependence: World Politics in Transitions}.
\newblock Little, Brown \& Co., Boston.

\bibitem[{Powers(1985)}]{Pow:85}
Powers, T. (1985).
\newblock Is there a way out?
\newblock \emph{Harpers}, 35--47.

\end{thebibliography}

\end{document}